\documentstyle{amsart}

\numberwithin{equation}{section}

\newtheorem{theorem}{Theorem}[section]

\newtheorem{lemma}[theorem]{Lemma}
\newtheorem{corollary}[theorem]{Corollary}
\newtheorem{proposition}[theorem]{Proposition}

\theoremstyle{remark}
\newtheorem{remark}[theorem]{Remark}

\newcommand{\cx}{{\Bbb C}}

\begin{document}

\title[Complexification and hypercomplexification of manifolds]{Complexification
and hypercomplexification of manifolds with a linear
connection}

\author{Roger Bielawski}
\thanks{Research supported by an EPSRC Advanced Research Fellowship}


\address{Department of Mathematics, University of Glasgow, Glasgow G12 8QW, UK }

\email{R.Bielawski@@maths.gla.ac.uk}

\begin{abstract} We give a simple interpretation of the adapted complex structure of
Lempert-Sz\H{o}ke and Guillemin-Stenzel: it is given by a polar decomposition of the complexified
manifold. We then give a twistorial construction of an $SO(3)$-invariant hypercomplex
structure on a neighbourhood of $X$ in $TTX$, where $X$ is a real-analytic manifold
equipped with a linear connection. We show that the Nahm equations arise naturally in this
context: for a connection with zero curvature and arbitrary torsion, the real sections
of the twistor space can be obtained by solving Nahm's equations in the Lie algebra
of certain vector fields. Finally, we show that, if we start with a metric connection,
then our construction yields an $SO(3)$-invariant hyperk\"ahler metric.

\end{abstract}

\maketitle

\thispagestyle{empty}

Let $X$ be a manifold equipped with a linear connection $\nabla$. Then the
tangent bundle $TX$ of $X$ has a canonical foliation (nonsingular on $TX\backslash X$)
by tangent bundles to
geodesics, i.e. by surfaces $T_\gamma \gamma$, where $\gamma$ is a
geodesic.
\par
A complex structure on $TX$ (or some neighbourhood of $X$) is called {\it
adapted} (to $\nabla$) if the leaves of the canonical foliation are complex
(immersed) submanifolds of $TX$.
\par
 For Riemannian connections the adapted complex structures were constructed and
studied by R. Sz\H{o}ke and L. Lempert \cite{Le,Sz,LSz}. An equivalent
definition was given by V. Guillemin and M. Stenzel \cite{GSt}. The results of
Lempert and Sz\H{o}ke can be formulated as follows:
\begin{theorem}\cite{Le,LSz,Sz,GSt} Let $(X,g)$ be a Riemannian manifold. There
exists a (unique) adapted complex structure for the Levi-Civita connection on
some neighbourhood of $X$ in $TX$ if and only if $(X,g)$ is a real-analytic
Riemannian manifold (i.e. both $X$ and $g$ are real-analytic).\end{theorem}
In this paper we shall give a simple construction of an adapted complex structure
which works for any real-analytic manifold $X$ with real-analytic linear
connection $\nabla$. It is sort of a polar decomposition of the complexified
manifold and generalises the basic example of the adapted complex structure on
$TG$, where $G$ is a compact Lie group equipped with the bi-invariant metric.
Actually, this construction is implicitly used by R. Sz\H{o}ke \cite{Sz} in the
proof that his (different) construction yields the adapted complex structure.
\par
In the Riemannian case, $(X,g)$, once we have the adapted complex structure $J$
on $TX$, we obtain a K\"ahler metric on $TX$ whose K\"ahler form is the canonical
$2$-form on the tangent bundle of a Riemannian manifold.  This metric
has in particular the property of being flat on leaves of the canonical
foliation.
\par
We would like to have an analogous statement in the non-Riemannian case,
i.e. obtain a K\"ahlerian (an almost complex connection whose curvature is
of type $(1,1)$) connection $\hat{\nabla}$ on $TX$, extending $\nabla$ and
such that the leaves of the canonical foliation are flat and totally
geodesic. This turns out to be a byproduct of our next construction: i.e.
hypercomplexification of $X$. Again, the motivating example is a compact
group with its bi-invariant metric: Kronheimer \cite{Kr} showed that there
exists a hyperk\"ahler metric on $T^{1,0}G^\cx$. We would like to construct
a hypercomplex structure on a neighbourhood of $X$ in $T^{1,0}TX$ where $X$
is any real-analytic manifold with a real-analytic linear connection
$\nabla$, and $TX$   is equipped with the adapted complex structure. We do
this in sections \ref{hyper} and \ref{real}. The Obata connection of this
hypercomplex structure restricted to $TX$ is the desired connection
$\hat{\nabla}$ (more precisely, we have to change the torsion if $\nabla$
was not torsion-free).
\par
If our original connection was (pseudo)-Riemannian, then we obtain a
(pseudo)-hyperk\"ahler metric a neighbourhood of $X$ in $T^{1,0}TX$. The restriction
of this metric to $TX$ is the K\"ahler metric of Lempert and Sz\H{o}ke.
\par
We should remark here that D. Kaledin \cite{Kal} and B. Feix \cite{Feix} showed
that, if  $M$ is a complex manifold with a real-analytic K\"ahlerian connection
$\nabla^\prime$, then there exists a canonical hypercomplex structure on a
neighbourhood of $M$ in $T^{0,1}M$. By the uniqueness result of Kaledin
\cite{Kal}, our construction gives the same hypercomplex structure as those of
Feix or Kaledin when applied to $M=TX$ with the adapted complex structure and
$\nabla^\prime=\hat{\nabla}$.

\par
 The hypercomplex and hyperk\"ahler structures obtained
here have a bigger symmetry than the more general ones of Feix and Kaledin.
 Namely, they admit an $SO(3)$ action rotating the complex structures.
\par
If the construction of the adapted complex structure should be viewed as a polar
decomposition of $X^\cx$, then the construction of the adapted hypercomplex
structure  should be viewed as solving certain Riemann-Hilbert problem. We make
this point of view precise in the case of connections with zero curvature and
arbitrary torsion. We show that the hypercomplex structure on $TTX$ can be in
this case obtained by solving a factorization problem in the loop group of the
group of holomorphic diffeomorphisms of $TX$. This yields a construction of this
hypercomplex structure through solving Nahm's equations in the Lie algebra of
vector fields on $TX$ invariant under the canonical involution $v\mapsto -v$ on
TX.

\section{Adapted complex structures\label{adapted}}

Let $X$ be a real-analytic manifold equipped with a  real-analytic linear
connection $\nabla$. We wish to construct an {\em adapted} complex structure on
a neighbourhood of $X$ in $TX$. Here ``adapted" means that the immersed surfaces
$T_\gamma \gamma$, where $\gamma$ is the  a
geodesic are complex submanifolds.
 \par
Our construction is sort of a polar decomposition of the complexification of
$X$.
\par
Let $X^{\Bbb C}$ be a complexification of $X$, i.e. a complex manifold
equipped with an antiholomorphic involution $\tau$ such that $\bigl(X^{\Bbb
C}\bigr)^\tau\simeq X$. Since $X$ is real-analytic, we can construct such a
complexification by a holomorphic extension of real-analytic transition
functions. As the connection $\nabla$ is real-analytic and $X$ is totally
real in $X^{\Bbb C}$, we can extend $\nabla$ to a holomorphic linear
connection $\nabla^{\Bbb C}$ on $T^{1,0}X^{\Bbb C}$ (i.e. a splitting of
$T^{1,0}\bigl(T^{1,0}X^{\Bbb C}\bigr)$ into vertical and horisontal holomorphic
bundles). The exponential map for this connection, which we
denote by $\exp^\cx$, is a holomorphic extension of the exponential map for
$\nabla$.
\par
At points $m$ of $X\subset X^{\Bbb C}$ we have the induced antiholomorphic
linear map $\tau_\ast:T_mX^{\Bbb C}\rightarrow T_mX^{\Bbb C}$. We denote by
the same symbol the corresponding linear map on $T_m^{1,0} X^{\Bbb C}$, i.e.
\begin{equation}\tau_\ast
(v-iIv)=\tau_\ast(v)-iI\tau_\ast(v),\label{tau}\end{equation} where $I$ is
the complex structure of $X^{\Bbb C}$ and $v\in TX^{\Bbb C}$ (real tangent
vector). Let $V^+_m$ and $V^-_m$ denote the $\pm 1$ eigenspaces of
$\tau_\ast$. They are interchanged by the multiplication by $i$. Thus we
have two bundles, $V^+$ and $V^-$, over $X$ which are subbundles of $T^{1,0}
X^{\Bbb C}$. Both $V^+$ and $V^-$ are canonically isomorphic to $TX$.
Moreover the restriction of $\exp^{\Bbb C}$ to $V^+$ is the exponential map
on $TX$.  On the other hand, as the differential of an exponential map is
identity at the zero section of a manifold,  $\exp^\cx$ restricted to $V^-$
is a diffeomorphism (near the zero section of $V^-$) between $V^-$ and
$X^\cx$. We can define a complex structure on $TX$ by pulling back the
complex structure of $X^{\Bbb C}$ to $V^-$ via this diffeomorphism. We have
\begin{proposition} The complex structure on $TX\simeq V^-$ obtained via the
(local) diffeomorphism
\begin{equation} \exp^\cx_{|_{V^-}}:V^-\rightarrow
X^\cx\label{iso}\end{equation} is a complex structure adapted to the
connection $\nabla$.
\end{proposition}
\begin{pf}
Since $\exp^\cx$ is a holomorphic extension of $\exp$, it sastisfies
\begin{equation}\exp^\cx\left((z_1+z_2)v\right)=
\exp^\cx\left(z_1\left(\frac{\partial}{\partial
z}\exp^\cx(zv)\right)_{|_{z=z_2}}\right)\label{group}\end{equation} for any
complex numbers $z_1$ and $z_2$ and for any $(1,0)$-tangent vector $v$ (as
long as the expressions are defined). Let us identify $V^+$ with $TX$, so
that $V^-_m=\{iv;v\in V^+_m\}$. Let $v$ be a tangent vector in $T_mX$ and
let $\exp(tv)$ be the corresponding geodesic. The image of the tangent
bundle to this geodesic under the map \eqref{iso} is the set $$
\exp^\cx\left(iq\frac{d\exp(tv)}{dt}_{|_{t=p}}\right)$$ where $p,q\in {\Bbb
R}$. As $\exp$ and $\exp^\cx$ coincide on $V^+$ and using \eqref{group},
this is the same as $\exp^\cx\left((p+iq)v\right)$. Hence the immersion sending ${\Bbb C}$ to the
tangent bundle of a geodesic is holomorphic.
 \end{pf}

We remark that this construction provides a natural interpretation of the
result in \cite{LSz} that global existence of the complex structure adapted
to a metric connection $\nabla$ implies that the curvature of $\nabla$ is
nonnegative. Indeed, global existence means that $\exp^\cx$ doesn't have
conjugate points along geodesics starting at $X$ and going in the {\em
imaginary} directions. Thus, in particular, we expect that the Jacobi
operator $R(\cdot,iv)iv$ is nonpositive at points of $X$.

\section{Adapted hypercomplex structure\label{hyper}}

We wish to construct a canonical hypercomplex structure on a neighbourhood
of $X$ in $TTX$. We shall do this by constructing a twistor space, i.e. a
complex manifold $Z$ fibering over ${\Bbb C}P^1$ equipped with a real
structure, i.e. an antiholomorphic involution $\sigma$ covering the antipodal
map on ${\Bbb C}P^1$, and a family of $\sigma$-invariant sections whose normal
bundle is a direct sum of ${\cal O}(1)$'s.
\par
Let $M$ be the complex manifold defined in the previous section, i.e. $M$ is
a neighbourhood of $X$ in $TX$ equipped with an adapted complex structure.
$M$ has  an antiholomorphic involution $\sigma=\tau^\ast$ such that
$M^\sigma=X$ and a holomorphic connection $\nabla^\cx$ which is an extension
of the original linear connection $\nabla$ on $X$. With these data we can
construct a twistor space $Z$. In this section we shall assume only that we
are given a complex connection on $M$; the real structure will be introduced
in the next section.

\par
$Z$ will be defined by gluing two open subsets of ${\Bbb C}\times T^{1,0}M$.
For any number $r>0$, let  $D_r=\{z\in {\Bbb C};|z|\leq r\}$ and $U_r$ be an
open subset of $T^{1,0}M$ defined as the set of these $v$ for which
$\exp^\cx (zv)$ is defined for all $z\in D_r$. Here $\exp^\cx$ denotes the
(holomorphic) exponential map for the connection $\nabla^\cx$.
\par
$Z$ is obtained by taking two copies of $U_2\times D_2$, parameterised by
$\{(\beta,y,\zeta); \beta\in T^{1,0}_yM, \zeta \in D_2\}$ and
$\{(\tilde{\beta},\tilde{y},\tilde{\zeta}); \tilde{\beta}\in
T^{1,0}_{\tilde{y}}M, \tilde{\zeta}\in D_2\}$, and identifying over
$1/2<|\zeta|<2$ via
\begin{equation}\tilde{\zeta}=\zeta^{-1}, \quad
\tilde{y}=\exp^\cx(-\beta/\zeta),
\quad \tilde{\beta}=\zeta^{-2}\frac{\partial}{\partial z}
\exp^\cx(z\beta)_{|_{z=-\zeta^{-1}}}.\label{transition}\end{equation}
We observe that the points of $M$ give rise to sections of $Z\rightarrow
{\Bbb C}P^1$: $m\mapsto (0,m,\zeta)$, $m\in M$. We have
\begin{proposition}
 The normal bundle  of sections corresponding to points $m$ of
 $M$ splits as the direct sum of line bundles of degree $1$.
\label{normal}\end{proposition}
\begin{pf}
Choose local holomorphic coordinates $z_1,\dots,z_n$ on $M$ and let
$z_1,\dots,z_n$,$z_{n+1}$,\linebreak $ \dots,z_{2n}$ be the induced
coordinates on $T^{1,0}M$, i.e. $z_{n+i}(Y)=dz_i(Y)$.  If $t\mapsto Y(t)$ is
an integral (complex) curve of the geodesic flow, then in the above local
coordinates $$ \dot{z}_k(t)=z_{n+k}(t),$$ $$ \dot{z}_{n+k}(t)=-\sum_{i,j}
\Gamma^k_{ij}(\pi\circ Y(t))z_{n+i}(t)z_{n+j}(t).$$ The transition functions
\eqref{transition} identify  $\tilde{z}_k$ with $z_k(-1/\zeta)$ and
$\tilde{z}_{n+k}$ with $\zeta^{-2}z_{n+k}(-1/\zeta)$. It follows that at a
section corresponding to a point of $m$, i.e. one where $z_{n+k}(t)\equiv
0$, $k=1,\dots,n$, the normal bundle has the transition function
\begin{equation} \left(\begin{matrix} 1 & -\zeta^{-1}\\ 0 &
\zeta^{-2}\end{matrix}\right),\label{trans}\end{equation} where the blocks
have size $n\times n$ (and correspond to the choice of coordinates
$z_1,\dots,z_n$ and $z_{n+1}, \dots,z_{2n}$). This means that the normal
bundle is isomorphic to $n$ copies of a rank $2$ bundle on $\cx P^1$ whose
transition function is \eqref{trans} with $(1\times 1)$-blocks. Therefore
this rank $2$ bundle $E$ is a nontrivial extension $$0\rightarrow {\cal
O}\rightarrow E\rightarrow {\cal O}(2)\rightarrow 0$$ and hence $E\simeq
{\cal O}(1)\oplus {\cal O}(1)$.
\end{pf}

In the above situation, i.e. when given a complex manifold $Z^{2n+1}$
fibering over ${\Bbb C}P^1$, it is well known that the parameter space $W$
of sections of $Z$ with normal bundle ${\cal O}(1)\otimes {\Bbb C}^{2n}$ is a
complex $4n$-dimensional manifold equipped with a canonical torsion-free
holomorphic connection known as the Obata connection (see, e.g., \cite{Feix2}).
Since, by the previous
lemma, $M$ is a complex submanifold of $W$, we can restrict the Obata
connection to $M$. The following comes as no surprise:
\begin{proposition}
 The restriction of the Obata connection to $M$ is the unique torsion-free
 connection with the same geodesics as the original connection $\nabla^{\Bbb
C}$.
\end{proposition}
\begin{pf}
It is enough to prove that every $\nabla^\cx$-geodesic in $M$ is a also a
geodesic in $W$ (equipped with the Obata connection).  Let $\gamma\subset M$
be a smooth geodesic disc (i.e. the image of a disc in ${\Bbb C}$ under a
$\nabla^\cx$-geodesic) containing a point $m\in M$. We can glue together two
copies of $T^{1,0}\gamma\times D_r$ using transition functions
\eqref{transition} and so obtain a twistor subspace $Z_\gamma$  of the
twistor space $Z$. Any section of $Z_\gamma$ is a section of $Z$. Moreover such
a section has the correct normal bundle in $Z$, at least in a neighbourhood of a
section $m(\zeta)$ corresponding to a point in $M$, since $W$ is a complete
family of sections at
$m(\zeta)$.   This means that the space of sections of $Z_\gamma$ is a
(complex) hypercomplex submanifold of $W$ and hence totally geodesic for the
Obata connection. However as $\gamma$ is a flat manifold (being a geodesic),
the hypercomplex structure of $T^{1,0}\gamma$ obtained from
\eqref{transition} is simply that of $\gamma\times \cx$ and hence $\gamma$
is a geodesic in $W$.
\end{pf}
\begin{remark} The proof shows that $M$ is totally geodesic in $W$. This follows
also from the existence of a ${\Bbb C}^\ast$ action on $Z$: $\lambda\cdot
(\beta, y, \zeta)=(\lambda\beta, y, \lambda\zeta)$ (more precisely the action may exist
only for $\log|\lambda|\leq c$, for some $c$), which induces an action
on the space of sections, i.e. $W$, whose fixed point set is $M$.
\end{remark}
\begin{remark}
The above proof also shows that the twistor space $Z$ has a property analogous
to the
one of adapted complex structure. Namely it is foliated by $3$-dimensional
complex manifolds $Z_\gamma$ fibering over ${\Bbb C}P^1$ and corresponding to
geodesics of $M$. These $Z_\gamma$ are "trivial" $3$-dimensional twistor spaces,
i.e. isomorphic to the twistor space of an open subset of $\cx^2$ with its
canonical hypercomplex structure. They give rise to flat $4$-dimensional
submanifolds of the space of sections $W$, but $W$ is not foliated by
these.\end{remark}

It is the last property that gives rise to a $PSl(2,{\Bbb C})$-action
on $W$.

\begin{proposition}
 The twistor space $Z$ admits a (local) action of $PSl(2,{\Bbb C})$ covering the
standard
 action on ${\Bbb C}P^1$. For each point $\zeta \in \cx P^1$, the subgroup
$\cx^\ast$
 fixing $\zeta$ acts on the fiber $Z_{\zeta}\simeq T^{1,0}M$  by multiplications
on the
 fibers of $T^{1,0}M$.
\end{proposition}
\begin{pf}
On the patch $\zeta\neq \infty$, the group of linear fractional
 transformations acts, where defined, by
\begin{equation} (m,\beta,\zeta)\mapsto
 \left(\exp^\cx_m\left(-\frac{c}{c\zeta +d}\beta\right),
\frac{\beta}{(c\zeta+d)^2},\frac{a\zeta+b}{c\zeta+d}\right).
\label{action0}\end{equation}
Here $1=ad-bc$ and this is helpful when checking that we do have a
$PSl(2,\cx)$ action on the first term. This action arises from expressing the
canonical (diagonal) action of $PSl(2,\cx)$ on ${\cal O}(1)\oplus {\cal O}(1)$
(the twistor space of a $Z_\gamma$ in Remark 2.4) as the lift of the canonical action
on ${\cal O}(2)$ in the extension
$$0\rightarrow {\cal
O}\rightarrow {\cal O}(1)\oplus {\cal O}(1) \rightarrow {\cal O}(2)\rightarrow 0$$
given by \eqref{trans}.
\end{pf}

\begin{corollary}
There exists a neighbourhood of $X$ in $W$ which admits an action of $SO(3,\Bbb R)$
induced from the local action of  $PSl(2,{\Bbb C})$ on $Z$.
\end{corollary}
\begin{pf}
From the previous proposition there exists a Lie subalgebra of vector fields on $W$
isomorphic to ${\frak so}(3,{\Bbb R})$. These vector fields vanish at points of $X$.
By the properties of ODE's we can find a neighbourhood of $X$ such that the integral
curves of vector fields corresponding to vectors in the unit sphere of
 ${\frak so}(3,{\Bbb R})$ exist for time $\pi$. We obtain action of $SU(2)$ with
 the center acting trivially because of \eqref{action0}.
\end{pf}

\section{ Real structures\label{real} }

We shall now assume that the complex manifold $M$ admits an anti-holomorphic
involution $\tau$ whose fixed points set is a manifold $X$ with $\dim_{\Bbb
R} X=\dim_\cx M$. We can extend $\tau_\ast$ to $T^{1,0}M$ by the formula
\eqref{tau}. Then $\tau_\ast$ is an antiholomorphic involution on the
complex manifold $T^{1,0}M$. We assume that the complex connection
$\nabla^\cx$ is compatible with $\tau_\ast$, which in particular means that
the exponential mapping $\exp^\cx$ commutes with $\tau_\ast$. These
conditions imply that $M$ is a neighbourhood of $X$ in $TX$ equipped with the
adapted complex structure for the restriction of $\nabla^\cx$ to $X$. In
particular, the geodesics of $\nabla^\cx$ passing through points of $X$ are
leaves of the canonical foliation of $TX$.
\par
Using $\tau_\ast$ we define a real structure on $Z$, i.e. an antiholomorphic
involution $\sigma$
covering the antipodal map on ${\Bbb C}P^1$. With the notation of
\eqref{transition},we put
\begin{equation}
\sigma(\beta,y,\zeta)=\left(-\tau_\ast(\tilde{\beta}),\tau_\ast(\tilde{y}),
-1/\bar{\zeta}\right).\label{sigma}\end{equation}

Proposition \ref{normal} implies that points of $X$ give rise to
$\sigma$-invariant sections of $Z$ with the normal bundle splitting into line
bundles with first Chern class $1$. Thus we obtain a hypercomplex structure on a
neighbourhood $W$ of $X$ in $T^{1,0}M$. We can summarize the properties of this
hypercomplex structure as follows:
\begin{proposition}  \begin{itemize}
\item[(1)] There exists an $SO(3,{\Bbb R})$-action on $W$ rotating the complex
structures. The fixed point set of this action is $X$.
\item[(2)] With respect to any complex structure, $W$ is biholomorphic to a
neighbourhood of $X$ in $T^{1,0}M$, where $M$ is a neighbourhood of $X$ in
$TX$ equipped with the adapted complex structure.
\item[(3)] For any leaf $L$ of the canonical foliation of $M$, its
holomorphic tangent bundle $T^{1,0}L$ is a hypercomplex submanifold of $W$, and
the hypercomplex structure of  $T^{1,0}L$ is the one of an open subset of ${\Bbb
C}^2$.
\item[(4)] $X$ is totally geodesic in $W$ and the Obata connection restricted to
$X$ is the unique torsion-free connection with the same geodesics as the
original connection $\nabla$ on $X$.\end{itemize}\label{properties}
\end{proposition}

We are not able to describe the hypercomplex structure on $T^{1,0}M$
directly in terms of the geometry of $X$ (apart from the case considered in
the next section). Such a description involves finding the full
$4n$-dimensional family of real sections. We can, however, describe the
$2n$-dimensional family corresponding to points of $M$ (where $M$ is a complex
submanifold of the fiber of $Z$ over $\zeta=0$). Indeed, as $M$ is
$TX$ with the adapted complex structure, the description of section
\ref{adapted} shows that any point of $m$ can be written uniquely as
$\exp^\cx_x(iv)$, where $x\in X$ and $v\in T_x X$. Define  tangent vectors
$V$ at $m$ and $\tilde{V}$ at $\tau(m)$ by
$V=\frac{d}{dt}\exp^\cx_x(itv)_{|_{t=1}}$,
$\tilde{V}=\frac{d}{dt}\exp^\cx_x(itv)_{|_{t=-1}}$. Then the following is a
real section of $Z$:
\begin{equation}\zeta\mapsto (2\zeta V,m,\zeta),\qquad
\tilde{\zeta}\mapsto \bigl(-2\tilde{\zeta}
\tilde{V},\tau(m),\tilde{\zeta}\bigr).\label{M}\end{equation}
\par

We now restrict the Obata connection to $M$, i.e. to the sections of the above
form. We have
\begin{proposition} The restriction of the Obata connection to $M$ is a
K\"ahlerian connection (i.e. an almost complex connection whose curvature is
a $(1,1)$-form) such that the leaves of the canonical foliation are flat and
totally geodesic.
\end{proposition}
\begin{pf}
The first part is obvious since the Obata connection is K\"ahlerian and the
second part follows from Proposition \ref{properties}(3).\end{pf}
\begin{remark} It follows now, from the uniqueness result of Kaledin, that the
hypercomplex structure on $T^{1,0}M$ arises also via the construction of Kaledin
\cite{Kal} or that of Feix \cite{Feix2}, if we start with the above connection
on $M$.\end{remark}

\begin{remark} Since we know the sections corresponding to points of $M$, we
can, in principle, compute the restriction of the Obata connection to $M$ by
the method of Merkulov \cite{Mer}. This involves computing, for a section of
the form \eqref{M}, its second formal neighbourhood $N^{(2)}$ in $Z$, and of
the (unique) splitting of the sequence $0\rightarrow N\rightarrow
N^{(2)}\rightarrow S^2(N)\rightarrow 0$, where $N$ is the normal bundle of
the section.
\end{remark}

\section{ Flat connections and Nahm's equations}

We shall describe explicitly sections of the twistor space and hence the
hypercomplex structure in the case of a linear connection with zero
curvature (and non-trivial torsion).  Again, let $M$ denote a neighbourhood
of $X$ in $TX$ equipped with the adapted complex structure and the
complexified connection $\nabla^\cx$. Let $\exp^\cx$ be the exponential map
for this connection. For a tangent vector $v\in T^{1,0}_mM $ we denote by
$\hat{v}$ the holomorphic vector field in a neighbourhood of $m$ obtained by
translating $v$ parallelwise along geodesics of $\nabla^\cx$ (as the
curvature is zero, this is the same as parallel translation along arbitrary
curves). For a holomorphic vector field $V$ on $M$, we denote by $e^V$ the
(local) diffeomorphism of $M$ defined by $e^V(m)=m(1)$, where $m(t)$ is the
curve with $m(0)=m$ and $\dot{m}(t)=V_{|_{m(t)}}$. Vanishing of the curvature of
 $\nabla^\cx$
implies the following key fact:
\begin{equation} \exp^\cx_p(\hat{v}_{|_p})=e^{\hat{v}}(p),
\label{identity}\end{equation}
where $v$ is a tangent vector in $T_p^{1,0} M$ and $p$ is in a simply-connected
neighbourhood of $x$.
\par
Because of this identity, the problem of finding real sections of $Z$ can be
solved via factorization in the (local) group of holomorphic diffeomorphisms
of $M$. The sections of the twistor space $Z$ will be parameterised by
points $x$ of $X$ and triples of vectors $V_1,V_2,V_3$ in $T_xX$ (usually
only triples close to zero). Let us denote by the same symbol $V_k$ the
corresponding vector in $T^{1,0}M$, i.e. $V_k-iIV_k$. The idea is that the
section of $Z$ will be of the form
\begin{equation}
 \zeta\mapsto \left(\beta(\zeta),g_+(1,\zeta)x,\zeta\right),\quad
 \tilde{ \zeta}\mapsto \left(\tilde{\beta}(\tilde{\zeta}),
 g_-(1,\tilde{\zeta})x,\tilde{\zeta}\right)\label{section2}
\end{equation}
where $g_\pm$ are certain holomorphic diffeomorphisms of $M$ and
$\beta(\zeta),\tilde{\beta}(\tilde{\zeta})$ are vector fields on $M$ defined
by parallel translation (with respect to the connection $\nabla^\cx$) of
vectors $$ (V_2+iV_3)+2iV_1\zeta+(V_2-iV_3)\zeta^2, \quad (V_2+iV_3) \tilde{
\zeta}^2+2iV_1 \tilde{ \zeta}+(V_2-iV_3)$$ which are in $T^{1,0}_xM$. The
diffeomorphisms  $g_\pm$ are defined as solutions of the following
Riemann-Hilbert problem:
\begin{equation} e^{-t\beta(\zeta)/\zeta}g_+(t,\zeta)=g_-(\tilde{\zeta},t)
\end{equation}
for $t\in [0,1]$ with $ g_+(\zeta,0)= g_-(\tilde{\zeta},0)=1$. This
equation, together with \eqref{identity}, shows that \eqref{section2} really
defines a section of $Z$. The $R$-matrix method,  cf. \cite{Int}, implies
that $g_\pm$ are obtained from vector field-valued solutions $B_0,B_1,B_2$
to the Nahm equations:
\begin{equation}
 \dot{B}_0=[B_0,B_1],\quad  \dot{B}_1=[B_0,B_2], \quad
 \dot{B}_2=[B_1,B_2]\label{Nahm}
\end{equation}
with initial conditions $B_0(0)=\widehat{V_2}+\widehat{iV_3}$,
$B_1(0)=2\widehat{iV_1}$, $B_2(0)=\widehat{V_2} -\widehat{iV_3}$. We have
\begin{equation}
g_+^{-1}\dot{g}_+=\frac{1}{2}B_1(t)+B_2(t)\zeta, \quad g_-^{-1}\dot{g}_-=-
\frac{1}{2}B_1(t)-B_0(t)\tilde{\zeta} ,\quad g_+(\zeta,0)=
g_-(\tilde{\zeta},0)=1.\label{initial}
\end{equation}
Since the antiholomorphic involution $\tau_\ast$ on $T^{1,0}M$ is compatible
with the vector field bracket, i.e.
$\tau_\ast[X,Y]=[\tau_\ast(X),\tau_\ast(Y)]$, the solutions $B_i$ to
\eqref{Nahm} (with given initial conditions) satisfy $\tau_\ast(B_0)=B_2$,
$\tau_\ast(B_1)=-B_1$. Consequently the sections of $Z$ obtained this way
are real, i.e. $\sigma$-invariant.
\par
Thus it remains to show that this construction works in our given
infinite-dimensional setting, i.e. that for $V_1,V_2,V_3$ in a neighbourhood
of $0$ in $T_xX$, there exists a solution to the equations \eqref{Nahm} and
\eqref{initial} on all of $t\in[0,1]$ (for the equation \eqref{initial} we
assume that $\frac{1}{2}\leq |\zeta|\leq 2$). Alternatively we can prove that
solutions exist for small $t$ for any $V_1,V_2,V_3$. The existence of solutions to \eqref{Nahm} follows from the Cauchy-Kovalevskaya theorem, since the Nahm equations for holomorphic vector fields become in local coordinates a Cauchy's problem for first-order PDE's.

It remains to show that there exists a solution to
\eqref{initial} for small $t$.
\begin{lemma} Let $M$ be a real (resp. complex) manifold and $V(t)$,
$t\in[0,\epsilon]$ be a one-parameter family of real (resp. holomorphic) vector
fields in a neighbourhood of $x\in M$. Then there exists a one-parameter family
$g(t)$, $t\in[0,\epsilon^\prime]\subset [0,\epsilon]$ of real (resp.
holomorphic) diffeomorphisms of a (smaller) neighbourhood of $x$ satisfying
$$g^{-1}\dot{g}=V(t), \qquad g(0)=1.$$
\end{lemma}
\begin{pf} We observe that solving the equation $\dot{h}h^{-1}=-V(t)$, $h(0)=1$
is simply equivalent to finding, for each point $p$ near $x$, a curve $p(t)$
with $p(0)=p$ and $\dot{p}=-V(t)$. Such a solution exists locally by the usual
theorems (e.g. Peano's) on existence and uniqueness of solutions to ODE's. Now
$g=h^{-1}$ is the desired family of diffeomorphisms. \end{pf}

Thus we can indeed find the full $4n$-dimensional family of sections of $Z$ by
solving the Nahm equations \eqref{Nahm} in the Lie algebra of holomorphic vector
fields on $M$.

\section{Hyperk\"ahler metrics}
In this section we consider the case of a Riemannian connection $\nabla$ on $X$.
We shall show that the Obata connection of section \ref{real} is given by a
hyperk\"ahler metric and that the restriction of this metric to $M\subset X^\cx$ is
the K\"ahler metric of Lempert and Sz\H{o}ke mentioned in the introduction.
\par
Let $X$ be a real-analytic manifold with a real-analytic metric $g$ and let $\nabla$ be
the Levi-Civita connection of $g$. We recall that the tangent bundle
$\pi:TX\rightarrow X$ of a Riemannian manifold has a canonical $1$-form $\Theta$
defined by
\begin{equation}
 \Theta(a)=g\bigl(z,d\pi(a)\bigr), \qquad a\in T_z(TX). \label{form}
\end{equation}
As Lempert and Sz\H{o}ke show in \cite{LSz}, if $TX$ is equipped with the
adapted complex structure, then there is a K\"ahler metric on $TX$ whose K\"ahler form is
$d\Theta$.
\par
We now wish to construct a hyperk\"ahler metric on a neighbourhood of $X$ in $TTX$.
 Let $M$ and $Z$ be as in sections \ref{hyper} and  \ref{real}.
We need to give a holomorphic ${\cal O}(2)$-valued $\sigma$-invariant
 $2$-form on fibers of the twistor space $Z$. We extend metric $g$ to a complex metric
 $\tilde{g}$ on a neighbourhood $M$ of $X$ in $X^\cx$. In the trivialisation
 \eqref{transition}
 we define a $1$-form on each fibre $T^{1,0}M$ of $Z$ by the formula
 \begin{equation}
 \tilde{\Theta}(a)=\tilde{g}\bigl(z,d\pi(a)\bigr), \qquad a\in T^{1,0}_z(T^{1,0}M).
 \label{Theta}
 \end{equation}
 This gives rise to a holomorphic ${\cal O}(2)$-valued and $\sigma$-invariant
  $1$-form on the fibers of $Z$ and hence we obtain a $2$-form
  $\Omega=d\tilde{\Theta}$. From this we obtain a triple
  $\omega_1,\omega_2,\omega_3$ of $2$-forms on $W$ such that
  \begin{equation}
 \Omega=(\omega_2+i\omega_3)+2i\omega_1\zeta+(\omega_2-i\omega_3)\zeta^2.
 \label{Omega}
\end{equation}
Moreover, if $J_1,J_2,J_3$ denotes the hypercomplex structure on $W$, then
 $\omega_1(\cdot,J_1\cdot)= \omega_2(\cdot,J_2\cdot)=\omega_3(\cdot,J_3\cdot)$
  and this is a non-degenerate (as $\omega_2$ is non-degenerate) symmetric tensor
  on $W$, i.e. a (possibly indefinite) hyperk\"ahler metric on $W$.
We wish to show that this hyperk\"ahler metric
restricted to $M$ is the metric of Lempert and Sz\H{o}ke. From \eqref{Omega},
 we need to compute \eqref{Theta} on sections of the form \eqref{M} and show that
the term linear in $\zeta$ is equal to \eqref{form} (up to the factor of $2i$ and
up to adding a closed form).
 Let $b$ be a $(1,0)$ tangent vector to $M$ at $m$. Then the corresponding variation
 of the  section \eqref{M} is given by $s(\zeta)=(b,2\zeta A,\zeta)$
  where $A$ is the resulting variation of $V$ and,
 for each $\zeta$, $d\pi(s(\zeta))=b$. Hence
 \begin{equation}
\tilde{\Theta}(s(\zeta))=\tilde{g}\bigl(2\zeta V, b\bigr). \label{comp}
 \end{equation}
 Now recall that $m$ is given as $\exp^\cx_x(iv)$ where $x\in X$ and $v\in T_xX$.
 Therefore $b=u(1)$, where $u(t)$ is the Jacobi vector field along the (complex)
 geodesic $\gamma(t)=\exp^\cx_x(itv)$. The vector $V$ is simply $\dot{\gamma}(1)$.
 Thus the $\zeta$-linear term in \eqref{comp} is simply
 $2\tilde{g}(\dot{\gamma}(1),u(1))$.
\par
Now recall  that a Jacobi vector field can be written (also in the complex case) as
$u(t)=\lambda \dot{\gamma}(t) + t\mu \dot{\gamma}(t) +U(t)$, where $U(t)$ is
orthogonal to $\dot{\gamma}(t)$ for any $t$. Therefore
\begin{multline}
\tilde{g}(\dot{\gamma}(1),u(1))=(\lambda+\mu)\tilde{g}(\dot{\gamma}(1),
\dot{\gamma}(1))= (\lambda+\mu)\tilde{g}(\dot{\gamma}(0),\dot{\gamma}(0))=\\
\tilde{g}(\dot{\gamma}(0),u(0))+ \tilde{g}(\dot{\gamma}(0),\mu\dot{\gamma}(1))=
\tilde{g}(\dot{\gamma}(0),u(0))+ \tilde{g}(\dot{\gamma}(0),\dot{u}(0)).
\end{multline}
Since $\gamma(0)=x\in X$, the last expression is simply
 \begin{equation} ig(v,u(0))+ig(v,\dot{u}(0)).\label{U}\end{equation}
 The Jacobi vector field $u$ was obtained by varying  geodesics of the
 form $\exp^\cx_x(itv)$ and hence we conclude that $u(0)=w$ is real and
 $\dot{u}(0)$ is purely imaginary. On the other hand, \eqref{Omega} implies that
  the exterior derivative of \eqref{U}
 is purely imaginary. Therefore the second term in \eqref{U} is a closed $1$-form
 and $\omega_1=d\Theta$.
 \par
 We finally compute the signature of our hyperk\"ahler metric. Observe first that the
  above construction  works perfectly well for a metric $g$ of any signature
  $(p,q)$.
 \begin{proposition}
If the metric $g$ on $X$ has signature $(p,q)$, then the $SO(3)$-invariant
hyperk\"ahler metric $G$ on $W$ has signature $(4p,4q)$.
 \end{proposition}
 \begin{pf} It is enough to compute the signature at points of $X\subset W$. As $X$ is
 $SO(3)$-invariant and totally real for any complex structure and the $2$-forms
 $\omega_i$ are non-degenerate, $T_xW=T_xX\oplus J_1T_xX \oplus J_2T_xX \oplus J_3T_xX$
 at points $x\in X$. As the hyperk\"ahler metric $G$ is $J_i$-invariant and its
 restriction to $T_xX$ is $g$, the signature of $G$ restricted to each direct summand
 is $(p,q)$. Finally, the summands are mutually orthogonal for $G$, as $T_xX$ is
isotropic for each $2$-form $\omega_i$.
 \end{pf}

{\bf Note.} After this paper has been completed I received from Robert Sz\H{o}ke a preprint ``Canonical complex structures assoviated to connections and complexifications of lie groups", in which he also introduces the construction given in section 1 of this paper.

\end{document}